\documentclass[11pt,bezier,amstex]{article}  
\usepackage[textsize=small]{todonotes}       

\usepackage{color}              
\usepackage{graphicx}
\usepackage[]{amsmath}
\usepackage{amssymb}
\usepackage{hyperref}
\usepackage{enumerate, amsfonts, amsthm, bbm}

\definecolor{MyDarkBlue}{rgb}{0,0.08,0.50}
\definecolor{BrickRed}{rgb}{0.65,0.08,0}

\hypersetup{
colorlinks=true,       
    linkcolor=MyDarkBlue,          
    citecolor=BrickRed,        
    filecolor=red,      
    urlcolor=cyan           
}

\newtheorem{Lemma}{Lemma}[section]

\newtheorem{Theorem}[Lemma]{Theorem}

\newtheorem{Condition}[Lemma]{Condition}

\newcommand{\prob}{\mathbb{P}}

\newcommand{\FF}{\mathcal{F}}

\newcommand{\Zbold}{{\mathbb{Z}}}

\newcommand{\expec}{\mathbb{E}}

\newcommand{\e}{{\mathrm e}}

\newcommand{\R}{{\mathbb R}}

\newcommand{\eqn}[1]{\begin{equation} #1 \end{equation}}
\newcommand{\eqan}[1]{\begin{align} #1 \end{align}}
\newcommand{\lbeq}[1]{\label{#1}}

\newcommand{\sss}{\scriptscriptstyle}

\newcommand{\cluster}{{\cal C}}

\newcommand {\vep}{\varepsilon}
\newcommand\1{\mathbbm{1}}
\newcommand{\indic}[1]{\1_{\{#1\}}}

\newcommand{\nn}{\nonumber}
\newcommand{\conn}{\longrightarrow}
\newcommand{\connn}[1]{\stackrel{#1}{\conn}}

\newcommand{\Zd}{{\mathbb Z}^d}
\newcommand{\sinfty}{{\sss \infty}}

\numberwithin{equation}{section}

\newcommand{\wh}{\widehat}

\newcommand{\ra}{\rightarrow}
\newcommand{\mP}{\mathbb{P}}
\newcommand{\mE}{\mathbb{E}}
\newcommand{\floor}[1]{\lfloor #1 \rfloor}
\newcommand{\Z}{\mathbb{Z}}
\newcommand{\re}{\mathbb{R}}
\newcommand{\mc}[1]{\mathcal{#1}}
\newcommand{\N}{\mathbb{N}}
\newcommand{\hlf}{\frac{1}{2}}
\newcommand{\blank}[1]{}

\begin{document}

\author{
Remco van der Hofstad\footnote{Department of Mathematics and
Computer Science, Eindhoven University of Technology, P.O.\ Box
513, 5600 MB Eindhoven, The Netherlands. E-mail {\tt
rhofstad@win.tue.nl}}
\and
Mark Holmes\footnote{Department of Statistics,
The University of Auckland, Private Bag 92019, Auckland 1142,
New Zealand. E-mail {\tt mholmes@stat.auckland.ac.nz}}
}

\title{The survival probability and $r$-point functions in high dimensions}

\maketitle

\begin{abstract}
In this paper we investigate the \emph{survival probability}, $\theta_n$,
in high-dimensional statistical physical models, where $\theta_n$ denotes the probability that the
model survives up to time $n$. We prove that if
the $r$-point functions scale to those of the canonical
measure of super-Brownian motion, and if a certain
self-repellence condition is satisfied, then $n\theta_n\ra 2/(AV)$, where
$A$ is the asymptotic expected number of particles alive at time $n$,
and $V$ is the vertex factor of the model. Our results apply to
spread-out lattice trees above 8 dimensions, spread-out oriented percolation
above $4+1$ dimensions, and the spread-out contact process above
$4+1$ dimensions. In the case of oriented percolation, this reproves
a result by the first author, den Hollander and Slade (that was
proved using heavy lace expansion arguments),
at the cost of losing explicit error estimates.
We further derive several consequences of our result
involving the scaling limit of the number of particles
alive at time proportional to $n$. Our proofs are
based on simple weak convergence arguments.
\end{abstract}

\section{Introduction and results}
\label{sec-setting}

A celebrated result by Kolmogorov \cite{Kolm38} states
that the probability $\theta_n$ that a
Galton-Watson branching process with offspring distribution
having mean 1 and variance $\gamma$, starting
from a single initial particle, survives until time $n$ satisfies
$n\theta_n\ra 2/\gamma$ as $n\ra \infty$
(see also~\cite[Theorem II.1.1.]{Per02}).
A related classical result by Yaglom \cite{Yagl47} states that
the population size $N_n$ at time $n$ is such that, conditionally
on survival up to time $n$, the
random variable $n^{-1}N_n$ converges weakly to a random variable $Y$
having an exponential distribution with mean $\gamma/2$.
Thus, the probability of survival up to time $n$ decays like
$1/n$, while on the event of survival, the number of alive
particles grows proportional to $n$. In this paper, we
study extensions of this result, and their ramifications,
to general spatial statistical mechanical models in sufficiently
high dimensions.

We next define the \emph{scaling limit} of the particle
numbers for critical Galton-Watson trees. The probability of a particle surviving
is rather small, and in the literature, two constructions have been
investigated to resolve this problem. The first construction
to deal with the vanishing survival probability is to start
with a large number of particles, i.e., take $N_0=\lceil nx\rceil$,
where $x>0$. In this case, at any time $t>0$, the number of particles
at time $0$ whose lineage survives until time $t$ has an approximate Poisson
distribution with parameter $2x/\gamma$. Then, the process
$(N_{tn}/n)_{t\geq 0}$ converges in distribution to
\emph{Feller's branching diffusion} \cite{Fell51},
which is the unique solution to a stochastic differential equation describing
a continuous-state branching process (see also \cite{Lamp67}
for related results). The second construction to deal with
the vanishing survival probability is to multiply the
measure by a factor of $n$, making sure that the measure of the
event of survival to time proportional to
$n$ converges to a finite and positive limit.
Then, the process $(N_{tn}/n)_{t\geq 0}$ converges in distribution,
where the notion of convergence in distribution is defined
in terms of convergence of integrals of bounded continuous functions having
support on paths that survive up to time $\vep>0$.
The resulting measure is a $\sigma$-finite measure rather than a probability measure,
and is called the \emph{canonical measure} of the branching process in
reference to canonical measures appearing in infinitely
divisible processes (see e.g.~\cite{Kall73}).
We can retrieve a probability measure by `conditioning' the measure
on surviving up to time $1$.

While the two constructions are quite different, they are
closely related. Indeed, in
the first construction (conditionally upon survival to time $1$) take any of the Poisson
$2x/\gamma$ initial particles whose lineage survives until time $1$.
Then the distribution of its rescaled numbers of descendants
is identical to that in the canonical measure
conditioned to survive up to time $1$.

The models we consider will be \emph{spatial}.
Embedding the branching process into $\Z^d$,
with the initial particle located at the origin, $0\in \Z^d$, and
where the offspring of any given particle are independently located at
neighbors of that particle in $\Z^d$, we obtain a \emph{branching random walk}.
Since multiple occupancy can occur, the state of this process at time $n$
is best described by a (random) \emph{measure}, where the measure of any subset of $\R^d$
is the number of particles of generation $n$ located in that set.  With
appropriate rescaling of space, time, mass (associated to each particle),
and of the underlying law, we obtain a sequence of finite (no longer
probability) measures $\mu_n$. Watanabe \cite{Wata68} shows that the measures
$\mu_n$ converge weakly to a measure $\N_0$
on the space of measure valued paths $(X_t)_{t\ge 0}$ that survive
for positive time, i.e.~$S\equiv\inf\{t>0\colon X_t(1)=0\}>0$
(where $X_t(f)\equiv\int f dX_t$).  The measure $\N_0$ is called
the \emph{canonical measure of super-Brownian motion} and is $\sigma$-finite,
with $\N_0(S>\vep)=2/\vep$ for every $\vep>0$.  The notion of weak
convergence is defined with respect to the finite measures
$\N_0^{\vep}(\cdot)\equiv\N_0(\cdot,S>\vep)$ (see
e.g.~\cite{HolPer07}), and in particular
$n\theta_{\floor{nt}}\ra \gamma^{-1}\N_0(S>t)=2/(\gamma t)$.
See \cite{Daws93,Per02} for detailed surveys of
super-processes and convergence towards them, and
\cite{Dynk94,Ethe00,LeGa99b} for introductions to
super-processes and continuous-state branching processes.

In this paper, we study extensions of these results in the context of
general spatial statistical mechanical models in sufficiently
high dimensions that converge (or are conjectured to converge)
to super-Brownian motion (SBM) in the sense of \emph{convergence
of $r$-point functions}. Convergence of $r$-point functions means
that the (rescaled) joint moments of particle numbers and locations converge (to those of SBM). The use of $r$-point
functions has a long history and tradition in statistical physics.
The main result of this paper is that convergence of
$r$-point functions, subject to two conditions that are
valid in all our examples, implies that the classical results
by Kolmogorov, Yaglom and (to some extent) Feller hold as well. As such,
our result confirms that convergence of
$r$-point functions is a relevant and important notion (see also \cite{HolPer07}).

Let us introduce the general setting that we investigate.
Let $\prob$ denote the probability measure describing the law of our model.
All our models have a notion of \emph{intrinsic distance},
in which $x\connn{n}y$ means that the shortest path between
$x$ and $y$ has length $n$. Let $\Z_+=\{0,1,2,\dots\}$ and $\R_+=[0,\infty)$.  Then for $\vec{x}\in \Zbold^{d(r-1)}$
and $\vec{n}\in \Z_+^{r-1}$ (or $\vec{n}\in \R^{r-1}_+$ for models where time is continuous), we let
    \eqn{
    \label{r-pt}
    t^{\sss(r)}_{\vec{n}}(\vec{x})=\prob(0\connn{n_i} x_i\forall i=1, \ldots, r-1)
    }
denote the $r$-point function in the model. Further, for
$\vec{k}=(k_1,\dots,k_{r-1})\in ([-\pi,\pi]^{d})^{r-1}$, we let
    \eqn{
    \label{r-ptk}
    \wh{t}^{\sss(r)}_{\vec{n}}(\vec{k})=\sum_{\vec{x}\in \Z^{d(r-1)}}\e^{i\vec{k} \cdot \vec{x}}t^{\sss(r)}_{\vec{n}}(\vec{x})
    }
denote its Fourier transform, and
    \eqn{
    \label{sp-def}
    \theta_n=\prob(\exists x\in \Zbold^d \colon 0\connn{n} x)
    }
the survival probability. Let $A_n=\{x\colon 0\connn{n} x\}$,
$N_n=\#\{x\colon 0\connn{n} x\}$, and
$S_n=\{N_n>0\}=\{A_n\ne \varnothing\}$, so that $\theta_n=\mP(S_n)$.
When the underlying model is defined
in discrete time, we define $n\vec{t}$ to be the vector
$(\floor{nt_1},\dots,\floor{nt_r})$.

In this paper, we investigate the asymptotics of the survival
probability, assuming the asymptotic behavior of the $r$-point functions. These results apply to
(a) lattice trees; (b) oriented percolation; and (c) the contact process,
all above their (model-dependent) upper critical dimension, where the general
philosophy in statistical physics suggests that these models behave like branching
random walk. In particular, when the allowed connections are sufficiently spread out,
e.g.~where all vertices within distance $L\gg 1$ of a vertex are considered
to be neighbors of that vertex, the following condition holds as
a theorem for each of these models, above their respective critical
dimensions:

\begin{Condition}[Convergence of the $r$-point functions]
\label{cond-rpoint}
(a) There exist constants $A,V>0$ all depending on $L$ such that for each $r\ge 2$ and $\vec{t}\in \R_+^{(r-1)}$,
    \begin{eqnarray}
   \frac{1}{A(VA^2n)^{r-2}}\wh{t}^{\sss(r)}_{n\vec{t}}(0)\ra
   \wh{M}^{\sss (r-1)}_{\vec{t}}(0), \quad \text{as }n\ra \infty,
    	\label{eq:rpointlimits0}
    \end{eqnarray}
where the quantities $\wh{M}^{\sss (r-1)}_{\vec{t}}(0)$
are the joint moments of the total mass at times
$t_1, \ldots, t_{r-1}$ of the canonical measure of SBM.
In particular, $\wh{M}^{\sss (r-1)}_{t\vec{1}_{r-1}}(\vec{0})=t^{r-2}2^{-(r-2)}(r-1)!$.\\
(b) There exist constants $A,V,v>0$ all depending on $L$ such that for each $r\ge 2$, $\vec{t}\in \R_+^{(r-1)}$, and $\vec{k}\in \re^{d(r-1)}$,
    \begin{eqnarray}
   \frac{1}{A(VA^2n)^{r-2}}\wh{t}^{\sss(r)}_{n\vec{t}}\left(\frac{\vec{k}}{\sqrt{vn}}\right)\ra \wh{M}^{\sss (r-1)}_{\vec{t}}(\vec{k}), \quad \text{as }n\ra \infty,
    	\label{eq:rpointlimits}
    \end{eqnarray}
where the quantities $\wh{M}^{\sss (r-1)}_{\vec{t}}(\vec{k})$ are the Fourier transforms of the moment measures of the canonical measure of SBM.
\end{Condition}

Condition \ref{cond-rpoint}(a) is the weaker of the above conditions, and can be rephrased as
    \eqn{
    \label{mass-conv}
    n\mE\Big[\prod_{i=1}^{r-1} \big(N_{t_in}/n\big)\Big]
    \ra A(VA^2)^{r-2}\wh{M}^{\sss (r-1)}_{\vec{t}}(0),
    }
where $\wh{M}^{\sss (r-1)}_{\vec{t}}(0)$ are the limits
of the joint moments of population sizes of critical branching processes with
variance one offspring distributions. Note that the convergence
in \eqref{mass-conv} makes no assumption on the \emph{spatial} locations
of the particles involved, however the evolution of $N_n$ is affected by spatial interaction present in our models. Condition \ref{cond-rpoint}(b), which contains (a), can be rephrased as
    \eqn{
    \label{eq:rpointsbm}
    \mE_{\mu_n}\left[\prod_{j=1}^{r-1}X^{\sss(n)}_{t_j}(\phi_{k_j})\right]\ra \mE_{\N_0}\left[\prod_{j=1}^{r-1}X_{ t_j}(\phi_{k_j})\right],
    }
where $\phi_{k_j}(x)=\e^{i k_j\cdot x}$ for $k_j \in \re^d$ and $x\in \Z^d$, and where
    \eqn{
    \label{Xn-def}
    X^{\sss(n)}_{t}(f)=\frac{1}{VA^2n}\sum_{x \in A_{nt}}
    f(x/\sqrt{vn}), \quad \text{and} \quad \mu_n(\cdot)=nVA\prob(\cdot).
    }
Thus, Condition \ref{cond-rpoint}(b) states that certain moment measures of the rescaled processes
under the measure $\mu_n$ converge to those of the canonical measure
of SBM. Condition \ref{cond-rpoint}(b) is the condition that is typically
proved in the literature.

Before stating our main result, we start by formulating two further
conditions. Let
    \eqn{
    \label{cluster-def}
    \cluster(0)=\{(x,n)\colon 0\connn{n} x\}
    }
denote the \emph{oriented cluster} of $0\in \Zbold^d$, i.e., all
vertices $x\in \Zbold^d$ to which $0$ is connected, and we let
$|\cluster(0)|$ denote its size. In continuous-time
models, we instead take
    \eqn{
    |\cluster(0)|=\int_{0}^{\infty} \#\{x\colon 0\connn{t} x\}dt.
    }
We make two central assumptions on our high-dimensional models:

\begin{Condition}[Cluster tail bound]
\label{cond-cluster-tail}
There exists a constant $C_{\sss \cluster}$
such that
    \eqn{
    \label{cluster-tail-bd}
    \prob(|\cluster(0)|\geq k)\leq C_{\sss \cluster}/\sqrt{k}.
    }
\end{Condition}

\begin{Condition}[Self-repellent survival property]
\label{cond-self-repel}
Let $\FF_m$ be the $\sigma$-field generated by
the vertices at distance at most $m$ from $0$, i.e.~by
$\{(x,n)\colon 0\connn{n} x, n\leq m\}$. Then
there exists a constant $C_{\theta}$ such that,
with $N_m$ equal to the number of $x$ with $0\connn{m} x$, almost surely for every stopping time $M\le n$,
    \eqn{
    \label{self-repel-sp}
    \prob(A_M \conn n\mid \FF_M)\leq C_{\theta} N_M \theta_{n-M}.
    }
\end{Condition}
The cluster tail condition follows from the literature for all models
under consideration. The self-repellent survival property in \eqref{self-repel-sp}
turns out to be easy to check, and we shall do this below.
The first of our main results is the following theorem:

\begin{Theorem}
\label{thm:sp_conv}
When Conditions \ref{cond-rpoint}(a), \ref{cond-cluster-tail} and \ref{cond-self-repel}
hold, as $n\ra \infty$,
   \eqn{
    n\theta_n\ra 2/(AV),
    }
and for each $t>0$,
    \eqn{
    \mu_n(X^{\sss(n)}_t(1)>0)\ra \N_0(X_{t}(1)>0)=2/t.
   }
Consequently, conditionally on $N_{tn}>0$,
the finite-dimensional distributions of
$(N_{sn})_{s\geq t}$ converge to those
of Feller's branching diffusion started from an exponential
random variable with mean $A^2Vt/2$.
\end{Theorem}

For oriented percolation, the result reproves a result from
\cite{HofHolSla07a, HofHolSla07b} (but without the error estimates)
in a relatively simple way. See also
\cite{KozNac08,KozNac11,Saka04} for related results on
survival probabilities. Our set-up is rather general, so
that in the future, it might be applicable to percolation
and lattice animals as well.

Theorem \ref{thm:sp_conv} is particularly important,
since the combination of the convergence of the $r$-point functions (as
formulated in Condition \ref{cond-rpoint}(b)) and Theorem \ref{thm:sp_conv} imply (see \cite{HolPer07}) that $\{\mu_n\}_{n \ge 1}$
converge in the sense of \emph{finite-dimensional distributions} to $\N_0$.  This is the second of our main results.

\begin{Theorem}
\label{thm:fdd-conv}
When Conditions \ref{cond-rpoint}(b), \ref{cond-cluster-tail}
and \ref{cond-self-repel} hold, the finite-dimensional distributions of the
process $(X^{\sss(n)}_{t})_{t>0}$ under $\mu_n$
converge to those of $(X_{t})_{t>0}$ under the measure $\N_0$.
\end{Theorem}

We now present some examples. All of the examples involve a function $D\colon \Zd \to [0,1]$, with
$\sum_{x \in \Zd} D(x)=1$
that obeys the properties of Assumption~D
in \cite[Section 1.2]{HofSla02} (whose precise form is not important for the
present paper), together with \cite[Equation~(1.2)]{HofSla03b}. This assumption
involves a parameter $L\in \N$, which serves to spread out the connections
and which will be taken to be large.

\paragraph{Spread-out oriented percolation above $4+1$ dimensions.}
The spread-out oriented bond percolation model is defined as follows.
Consider the graph with vertices $\Zd \times \Z_+$ and with
directed bonds $((x,n),(y,n+1))$, for $n \in \Z_+$ and $x,y \in
\Zd$.  Let $p \in [0,\|D\|_{\sinfty}^{-1}]$, where $\|\cdot\|_{\sinfty}$ denotes
the supremum norm, so that $pD(x) \leq 1$ for all $x\in\Zd$. We
associate to each directed bond $((x,n),(y,n+1))$ an independent
random variable taking the value $1$ with probability $pD(y-x)$ and
the value $0$ with probability $1-pD(y-x)$. We say that a bond is
{\em occupied}\/ when the corresponding random variable is $1$ and
{\em vacant}\/ when it is $0$. The joint probability distribution of the bond variables
will be denoted by $\prob_p$, and the corresponding expectation by
$\mE_{\,p}$.

We say that $(x,n)$ is {\em connected}\/ to $(y,m)$, and write $(x,n)
\conn (y,m)$, if there is an oriented path from $(x,n)$ to $(y,m)$
consisting of occupied bonds. Note that this is only possible when
$m \geq n$. By convention, $(x,n)$ is connected to itself. We write
$(x,n) \conn m$ if $m \geq n$ and there is a $y \in \Zd$ such that
$(x,n) \conn (y,m)$. The event $\{(0,0) \conn \infty\}$ is the event that $\{(0,0)\conn n\}$
occurs for all $n$.  There is a critical threshold $p_c > 0$ such that
the event $\{(0,0) \conn \infty\}$ has probability zero for $p < p_c$
and has positive probability for $p>p_c$.
The survival probability at time $n$ is defined by
    \eqn{
    \label{thetadef}
    \theta_n(p) =\prob_p((0,0)\conn n),
    }
    and we let $\theta_n=\theta_n(p_c)$.
General results of \cite{BG90,GH02} imply that $\lim_{n \to \infty}
\theta_n=0$.

Then, for $\mP=\mP_{p_c}$,
Condition \ref{cond-rpoint} is proved in \cite{HofSla03b}.
Condition \ref{cond-cluster-tail} holds by \cite{BarAiz91, HofSla03b, NguYan93,NguYan95},
while Condition \ref{cond-self-repel}
follows from a union bound (i.e.~$\mP(\cup_{x\in A_M}\{x\ra n\}|\mc{F}_M)\le \sum_{x\in A_M}\mP(x\ra n|\mc{F}_M)$) and the strong Markov property.

\paragraph{Spread-out contact process above $4+1$ dimensions.}
We define the spread-out contact process as follows.  Let $\cluster_n\subset\Zd$
be the set of infected individuals at time $n\in\R_+$, and let
$\cluster_0=\{0\}$.  An infected site $x$ recovers in a small time interval
$[n,n+\vep]$ with probability $\vep+o(\vep)$ independently of $n$, where
$o(\vep)$ is a function that satisfies $\lim_{\vep\to0}o(\vep)/\vep=0$.
In other words, $x\in\cluster_n$ recovers at rate 1.  A healthy site $x$ gets
infected, depending on the status of its neighbors, at rate
$\lambda\sum_{y\in\cluster_n}D(x-y)$, where $\lambda\geq0$ is the infection rate.
We denote by $\prob^\lambda$ the associated probability measure.

By an extension of the results in \cite{BG90,GH02} to the spread-out
contact process, there exists a unique critical value $\lambda_c\in(0,\infty)$
such that
\begin{align}\lbeq{lambc-def}
\theta(\lambda)&\equiv\lim_{n\ra\infty}\prob^\lambda(\cluster_n\neq\varnothing)
 \begin{cases}
 =0,&\text{if }\lambda\leq\lambda_c,\\ >0,&\text{if }\lambda>\lambda_c,
 \end{cases}
\end{align}
and we define
    \eqn{
    \theta_n=\theta_n(\lambda_c)=\prob^{\lambda_c}(\cluster_n\neq\varnothing).
    }

Condition \ref{cond-rpoint} is proved in \cite{HofSak04, HofSak10}.
Condition \ref{cond-cluster-tail} holds
by \cite{HofSak04, HofSak10, BarWu98, Saka01}, while
Condition \ref{cond-self-repel}
again follows from a union bound and the strong Markov property.

\paragraph{Spread-out lattice trees above $8$ dimensions.}
 A lattice tree is a finite connected set of lattice bonds (and their associated end vertices) containing no cycles.  For fixed $z>0$, every such tree $T\ni 0$ with bond set $B$ is assigned a weight $W_{z}(T)=z^{|B|}\prod_{(x,y)\in B}D(y-x)$, and we define $\rho_z(x)=\sum_{T\ni 0,x}W_z(T)$.  The radius of convergence $z_c$ of $\sum_{x\in \Z^d}\rho_z(x)$ is finite.  Let $W(\cdot)=W_{z_c}(\cdot)$ and $\rho=\rho_{z_c}(0)$.  We define a probability measure on the (countable) set of lattice trees containing the origin by $\mP(T)=\frac{W(T)}{\rho}$.  Given a lattice tree $T\ni 0$, we define $A_n(T)=\{a_1,\dots, a_{N_n}\}$ to be the (ordered) set of vertices in $T$ of tree distance $n\in \Z_+$ from the origin under some arbitrary but fixed ordering of $\Z^d$.

Condition \ref{cond-rpoint} is the main result in \cite{Holm08b}.
Condition \ref{cond-cluster-tail} follows from the detailed
asymptotics for $\mP(|T|=n)\sim c n^{-3/2}$ proved in \cite{DerSla97, DerSla98}.
We next check Condition \ref{cond-self-repel}, for which it is enough to show that the result holds a.s.~for every deterministic time $m\le n$.
Letting $T_m$ denote the tree up to tree distance
$m$ from the root, we have that $\mP(A_m\conn n\mid T_m=\tau_m)$ is equal to
    \[
    \frac{W(\tau_m)}{\sum_{T\colon T_m=\tau_m}W(T)}\sum_{R_1\ni a_1}
    \dots \sum_{R_{N_m}\ni a_{N_m}}\prod_{i=1}^{N_m}W(R_i)\indic{R_i \text{ avoid each other and }\tau_m}\indic{\cup_j R_j \text{ survives at least until }n-m},
    \]
where $\sum_{ R\ni a}$ is a sum over lattice trees $R$ containing $a\in \Z^d$, and
we recall that $A_m=\{a_1, \ldots, a_{N_m}\}$.

The final indicator function is bounded above by
$\sum_{j}\indic{\mc{S}_{R_j}\ge n-m}$, where $\mc{S}_{T}$
is the survival time of $T$.  By taking the sum over $j$
outside and dropping the restriction that $R_j$ avoids other
$R_i$ and $\tau_m$, this is bounded above by
    \begin{eqnarray}
    &&\sum_{j=1}^{N_m}\sum_{R_j\ni a_j}W(R_j)I_{\{\mc{S}_{R_j}\ge n-m\}}\Bigg[\frac{W(\tau_m)}{\sum_{T\colon T_m=\tau_m}W(T)}\\
    &&\times \sum_{R_1\ni a_1}
    \dots \sum_{R_{j-1}\ni a_{j-1}}\sum_{R_{j+1}\ni a_{j+1}}\dots\sum_{R_{N_m}\ni a_{N_m}}
    \prod_{i\ne j}W(R_i)\indic{R_i, i\ne j \text{ avoid each other and }\tau_m}\Bigg]\nn\\
    &&\quad \le \sum_{j=1}^{N_m}\sum_{R_j\ni a_j}W(R_j)\indic{\mc{S}_{R_j}\ge n-m}=N_m\rho\theta_{n-m},\nn
    \end{eqnarray}
where we have used the fact that the interaction term makes the graph
$\tau_m \cup_{i\ne j}R_{i}$ a lattice tree $T$ with $T_m=\tau_m$, and
weight $W(T)=W(\tau_m) \prod_{i\ne j}W(R_i)$, so the numerator in brackets
is no more than the denominator.  This verifies Condition \ref{cond-self-repel}.

Our main results can be restated in terms of the above models as follows:

\begin{Theorem}
\label{thm-sp}
Let $L\gg 1$, and let $d>4$ for spread-out oriented percolation and the spread-out contact
process, and $d>8$ for spread-out lattice trees.  Then, with
$A,V,v>0$ all depending on $L$ such that for each $\vec{t}\in \R_+^{(r-1)}$ and $\vec{k}\in \re^{(r-1)}$

\begin{eqnarray}
\frac{1}{A(VA^2n)^{r-2}}\wh{t}^{\sss(r)}_{n\vec{t}}\Big(\vec{k}/\sqrt{vn}\Big)\ra \wh{M}^{\sss (r-1)}_{\vec{t}}(\vec{k}),\quad \text{as }n\ra \infty,
	\label{eq:rpointlimits2}
\end{eqnarray}
the asymptotics
    \eqn{
    n\theta_n\ra 2/(AV) \quad \text{and}\quad
    \mu_n(X^{\sss(n)}_t(1)>0)\ra \N_0(X_{t}(1)>0)=2/t, \quad \text{as }n\ra \infty,
    }
hold. As a consequence, the finite-dimensional distributions of the
process $(X^{\sss(n)}_{t})_{t>0}$ under $\mu_n$
converge to those of $(X_{t})_{t>0}$ under the measure $\N_0$.
\end{Theorem}

We close this section with two possible extensions to our results.

\paragraph{Long-range models.}
In all our models, we assume that $D$ has finite spatial variance, so
that SBM can arise as the scaling limit. In the literature,
\emph{long-range} models have attracted considerable attention.
See \cite{CheSak08, CheSak09, CheSak11} for results on long-range
oriented percolation, \cite{Heyd11} for long-range self-avoiding walk,
and \cite{HeyHofSak08} for percolation, self-avoiding walk and
the Ising model. In long-range models, the random walk step distribution $D$
has infinite variance. The simplest example arises when
    \begin{equation}
    \label{eqDefD}
    D(x)=\frac{(1+|x|/L)^{-(d+\alpha)}}{\sum_{y\in\Zd}(1+|y|/L)^{-(d+\alpha)}},
    \qquad x\in\Zd,
    \end{equation}
where $\alpha\in (0,2)$, and $|x|$ denotes the Euclidean norm of
$x\in \Zd$. The results in \cite{CheSak08, CheSak09, CheSak11}
suggest that the upper critical dimension of oriented percolation
equals $2\alpha$, while \cite{HeyHofSak08} indicates that it is
$3\alpha$ for percolation, and $2\alpha$ for self-avoiding walk and
the Ising model.

We believe that Condition \ref{cond-rpoint}(a)
holds for these models above their respective upper critical dimensions.
Once this is proved, Theorem \ref{thm:sp_conv} then implies convergence
of the survival probability in each case.
However, random walk with step distribution $D$ converges to
$\alpha$-stable motion rather than Brownian motion, a fact
that is proved to hold for self-avoiding walk above $2\alpha$
dimensions in \cite{Heyd11}. Therefore, Condition \ref{cond-rpoint}(b)
does \emph{not} hold, and should be replaced with
convergence towards the canonical measure of super-stable motion.

By considering branching random walks, where the population size process is independent of the random walk  step-distribution, it is easy to see that the law of the total mass process under the canonical measure of super-stable motion is the same as under $\N_0$.  Thus by \cite[Theorem 2.6]{HolPer07}, in the long range setting, convergence of the $r$-point functions and the survival probability still implies convergence in the sense of finite-dimensional distributions.  Therefore to prove a version of Theorem \ref{thm:fdd-conv}
in the long-range setting, it is sufficient to prove the convergence
of the $r$-point functions in Condition \ref{cond-rpoint}(b).

\paragraph{Spread-out percolation above $6$ dimensions.}
Let $p \in [0, \| D \|_\infty^{-1}]$ be a parameter.
We declare a bond $\{u,v\}$ to be {\em occupied}
with probability $p D(v - u)$ and {\em vacant} with probability
$1 - p D(v - u)$. The occupation status of all bonds are
independent random variables.
The law of the configuration of occupied bonds (at the critical
percolation threshold) is denoted by $\prob_{p_c}$ with
corresponding expectation denoted by $\expec_{p_c}$.
Given a configuration we say that $x$ is connected to $y$,
and write $x \connn{n} y$, if there is a path of occupied bonds from $x$ to $y$,
and the path with minimal number of bonds connecting $x$ and $y$ has
precisely $n$ edges. For percolation, Condition \ref{cond-rpoint} is not
known. The bound $\theta_n\leq C/n$
is proved in \cite{KozNac08}
(in fact, we use Condition \ref{cond-self-repel} together
with an adaptation of the argument in \cite{KozNac08}
to prove that $\theta_n\leq C/n$ in our general setting).
Condition \ref{cond-cluster-tail} follows from \cite{HarSla90a}
together with \cite{BarAiz91}, see also \cite{HarSla00a, HarSla00b}.
As a result, for percolation, our results hold as soon as
Condition \ref{cond-rpoint} is proved.

The above discussion suggests the following research program to
identify the right constants in arm-probabilities in high-dimensional
percolation, both in the intrinsic as well as in the
Euclidean or extrinsic distance: (1) prove the convergence
of the $r$-point functions in Condition \ref{cond-rpoint}(b)
(from which the right constant in the survival probability or
intrinsic one-arm probability would follow,
improving upon the results in \cite{KozNac08});
(2) prove \emph{tightness} for convergence towards SBM;
(3) identify the right constant for the \emph{extrinsic}
one-arm probability, improving upon the result in
\cite{KozNac11}. For the last step, an important
ingredient showing that it is unlikely
that a short path exists to the boundary of a Euclidean
ball is proved in \cite[Theorem 1.5]{HofSap11}.

\bigskip
The remainder of this paper is organized as follows. In Section \ref{sec-weak-bounds},
we prove an upper bound on $\theta_n$ that is of the correct order, but with the wrong constant.
In Section \ref{sec-Pf-main-thm},
we use weak-convergence arguments to identify the correct constant,
and prove the consequences of convergence of the survival probability.

\section{Weak upper bound on the survival probability}
\label{sec-weak-bounds}
The following theorem gives a weak upper bound on the survival probability.

\begin{Theorem}
\label{thm-weak-ub}
When Conditions \ref{cond-cluster-tail} and \ref{cond-self-repel}
 hold,
there exists a constant $c_+$ such that
    \eqn{
    \theta_n \leq c_+/n.
    }
\end{Theorem}

\proof We follow \cite{KozNac08}, where a similar bound was proved for
the intrinsic one-arm in percolation. We split $\theta_{4n}$ into two parts,
    \eqn{
    \label{split-theta}
    \theta_{4n}=\prob(N_m\geq \vep n~\forall m\in [n, 3n], 0\conn 4n)
    +\prob(\exists m\in [n, 3n]\colon N_m<\vep n, 0\conn 4n).
    }
We can bound the first probability using \eqref{cluster-tail-bd}, since
$|\cluster(0)|\geq 2\vep n^2$ if $N_m\geq \vep n$ for all $m\in [n, 3n]$. Therefore,
    \eqn{
    \prob(N_m\geq \vep n~\forall m\in [n, 3n], 0\conn 4n)
    \leq \prob(|\cluster(0)|\geq 2\vep n^2) \leq  \frac{C_{\sss \cluster}}{n\sqrt{2\vep}}.
    }
In the second probability in \eqref{split-theta}, we let $J\geq n$ be the first
$m\in [n, 3n]$ such that $0<N_m<\vep n$, and we condition on
$\mc{F}_J=\sigma((A_m)_{m\leq J})$. Then, by \eqref{self-repel-sp},
    \eqn{
    \prob(A_J\conn 4n\mid \mc{F}_J)
    \leq N_J C_{\theta} \theta_{n}\leq \vep n C_{\theta} \theta_{n}.
    }
As a result,
    \eqn{
    \prob(\exists m\in [n, 3n]\colon N_m<\vep n, 0\conn 4n)
    =\expec[\indic{n\le J\leq 3n}\prob(A_J\conn 4n\mid \mc{F}_J)]
    \leq \vep C_{\theta} n \theta_{n}^2,
    }
where we use the fact that $n\le J$ implies that $0\conn n$.
Thus, we end up with the inequality
    \eqn{
    \theta_{4n}\leq \frac{C_{\sss \cluster}}{n\sqrt{2\vep}}+\vep C_{\theta} n \theta_{n}^2.
    }
Take $\vep=c_2^{-4/3}$ and take $c_2>1$ so large that
    \eqn{
    \label{cond-constant}
   2^{-\hlf} C_{\sss \cluster}c_2^{2/3}+C_{\theta} c_2^{2/3}\leq c_2/4.
    }
Then, it is easy to prove by induction that $\theta_{4^k}\leq c_2 4^{-k}$ for every $k\ge 1$.
By monotonicity of $n\mapsto \theta_n$, this immediately implies that $\theta_n\leq (4c_2)/n$.
This completes the proof of Theorem \ref{thm-weak-ub}.
\qed

\section{Identifying the constant: Proof of Theorem \ref{thm:sp_conv}}
\label{sec-Pf-main-thm}
In this section, we make use of general weak convergence arguments to
prove that $n\theta_n\ra 2/(AV)$. We rely crucially on a result that is essentially a special case of \cite[Proposition 2.3]{HolPer07},
which requires the introduction of some more notation.  Let $M_F(\re^d)$ (resp.~$M_1(\re^d)$) denote the space of finite (resp.~probability) measures on $\re^d$ equipped with the topology of weak convergence.
Let $\mc{D}_G$ denote the set of discontinuities of a function $G$, and $D(E)$ denote the space of c\`adl\`ag $E$-valued functions with the Skorohod topology.  When we say that $\mu$ is a measure on (a topological space) $E$, this means that it is a measure with respect to the Borel $\sigma$-algebra on $E$.


\begin{Lemma}
\label{lem:HPmass}
Suppose that Condition \ref{cond-rpoint}(a) holds.
Then for every  $s,t,\eta>0$,
and every bounded Borel measurable
$H\colon \re\ra \re$ such that $\N_0(X_t(1)\in \mc{D}_H)=0$,
    \eqn{
    \label{eq:rpoint2}
    \expec_{\mu_n}\left[\indic{X_{s}^{\sss(n)}(1)>\eta}H(X^{\sss(n)}_t(1))\right]\ra
    \expec_{\N_0}\left[\indic{X_{s}(1)>\eta}H(X_t(1))\right],
    \quad \text{as }n\ra \infty.
    }
\end{Lemma}
\proof We follow the proof of \cite[Proposition 2.3]{HolPer07}.   For convenience, we drop the superscripts $~^{\sss (n)}$.

By Condition \ref{cond-rpoint}(a), $\{\mu_n\}_{n \ge 1}$ is a sequence of finite measures on
$D(M_F(\re^d))$ such that for every $r \ge 1$ and $\vec{t}\in [0,\infty)^r$, (\ref{eq:rpointsbm}) holds when $\phi_{k_j}=1$ for each $j$.

Fix $s,t,\eta>0$. Let $Y_s=X_s(1)$, and define $P_n=P_{n,s,t}\in M_1(\re^2)$ and $P=P_{s,t}\in M_1(\re^2)$ by
    \[
    P_{n}(A)=\frac{\expec_{\mu_n}[Y_s\indic{(Y_s,Y_t)\in A}]}{\expec_{\mu_n}[Y_s]},
    \quad \text{and } P(A)=\frac{\expec_{\N_0}[Y_s\indic{(Y_s,Y_t)\in A}]}{\expec_{\N_0}[Y_s]},
    \]
where these measures are well defined since
    \[
    \expec_{\mu_n}[Y_s]\ra \expec_{\mathbb{N}_0}[Y_s]\in (0,\infty).
    \]
On each of these spaces let $(W,Z)$ be the canonical random vector, i.e.~$(W,Z)(\omega_1,\omega_2)=(\omega_1,\omega_2)$.  Then, for every $m_1,m_2\ge 0$,
    \eqn{
    \label{eq:Pmoments}
    \expec_{P_n}\left[W^{m_1}Z^{m_2}\right]
    =\frac{\expec_{\mu_n}\left[Y_s^{m_1+1}Y_t^{m_2}\right]}{\expec_{\mu_n}[Y_s]}\ra \frac{\expec_{\N_0}\left[Y_s^{m_1+1}Y_t^{m_2}\right]}{\expec_{\N_0}[Y_s]}
    =\expec_{P}\left[W^{m_1}Z^{m_2}\right],
    }
i.e.~the moments of $(W,Z)$ under $P_n$ converge to those under $P$.

Furthermore (see e.g.~\cite[Lemma 4.1]{HolPer07}) there exists $\delta>0$ such that,
    \eqn{
    \expec_{P}\left[\e^{\delta (W+Z)}\right]=\frac{\expec_{\N_0}\left[Y_s\, \e^{\delta (Y_{s}+ Y_t)}\right]}{\expec_{\N_0}[Y_s]}<\infty,
    }
i.e.~the moment generating function of $(W,Z)$ under $P$ is finite
in a neighborhood of $(0,0)$.
It then follows (see e.g.~\cite[Theorems 30.1 and 30.2, and Problems 30.5 and 30.6]{Bill86}) that $P_n$ converges weakly to $P$, and therefore for $G\colon \re^2\ra \re$ bounded and such that $P((W,Z)\in \mc{D}_G)=0$,
    \[
    \expec_{P_n}[G(W,Z)]\ra \expec_{P}[G(W,Z)].
    \]
In other words, for each bounded $G\colon \re^2\ra \re$ such that $\N_0((Y_{s},Y_t)\in \mc{D}_G)=0$,
    \[
    \expec_{\mu_n}\left[Y_s G(Y_s,Y_t)\right]\ra \expec_{\N_0}\left[Y_s G(Y_s,Y_t)\right].
    \]
Let $H$ be as in the statement of the lemma, and define
    \[
    G_H(x,y)=\begin{cases}
    \frac{H(y)}{x}, & \text{ if }x>\eta\\
    0, & \text{otherwise}.
    \end{cases}
    \]
Then $G_H$ is bounded, and $\mc{D}_{G_H}=\{(x,y)\colon y\in \mc{D}_H \text{ or } x=\eta\}$, whence $\N_0((X_s,X_t)\in \mc{D}_{G_H})=0$.  The claim follows since $Y_s
G_H(Y_s,Y_t)=\indic{Y_s>\eta}H(Y_t)$.
\qed

\bigskip

\noindent {\em Proof of Theorem \ref{thm:sp_conv}.}  By Theorem \ref{thm-weak-ub}, we have that $n\theta_n$ is bounded. In order to investigate the limit of $n\theta_n$, we split, for each fixed $\vep>0$,
    \begin{eqnarray}
    \label{eq:fish}
    n\theta_n&=&n\mP(N_n>\vep n)+n\mP(0<N_n\le \vep n).
    \end{eqnarray}
The first term is equal to $(AV)^{-1}\mu_n(X_1^{\sss(n)}>c\vep)$, with $c=(VA^2)^{-1}$.
From Lemma \ref{lem:HPmass} with $s=1$, $\eta=c\vep$ and with the continuous function $H\equiv 1$
(and Condition \ref{cond-rpoint}(a)), we have that the first term on the right converges to
$(AV)^{-1}\N_0(X_1(1)>c\vep)$, and this converges to  $(AV)^{-1}\N_0(X_{1}(1)>0)=2/(AV)$ as $\vep\ra 0$.
Since $n\mP(0<N_n\le \vep n)\geq 0$, this immediately proves that
    \eqn{
    \label{weak-lb}
    \liminf_{n\rightarrow \infty} n\theta_n\geq 2/(AV).
    }

In order to identify the limit, we proceed as in \cite{HHS_IIC}.  Let $\delta\in (0,1)$ and let $\{n_k\}=\{n_k(\delta)\}$ be any subsequence of $\N$ such
that $n_k\theta_{n_k}\ra \limsup_{n} n\theta_n=\overline{b}$, and $(1-\delta)n_k\theta_{(1-\delta)n_k}\ra b_{\delta}$ for some $b_{\delta}\ge 2/AV$.
This can be achieved by first taking
a subsequence $\{m_l\}$ for which $m_l\theta_{m_l}\ra \overline{b}$,
and then taking a further subsequence $\{m_{l_k}\}$ such that  $(1-\delta)m_{l_k}\theta_{(1-\delta)m_{l_k}}\ra b_{\delta}$.
The required sequence is then $n_k=m_{l_k}$.

Similarly to (\ref{eq:fish}), for $\delta, \vep, \vep'\in (0,1)$ we write
 \begin{eqnarray}
    \label{eq:fish2}
    n_k\theta_{n_k}&=&n_k\mP(N_{(1-\delta)n_k}>\vep n_k,N_{n_k}>\vep' n_k )\nn\\
    &&+n_k\mP(N_{(1-\delta)n_k}>\vep n_k,0<N_{n_k}\le \vep' n_k )+n_k\mP(0<N_{(1-\delta)n_k}\le \vep n_k,N_{n_k}>0)\nn\\
    &=&A_{k,\delta,\vep,\vep'} + B_{k,\delta,\vep,\vep'} + D_{k,\delta,\vep}.
    \end{eqnarray}
Since the above is true for each $\delta,\vep,\vep'$, it follows that also
   \begin{eqnarray}
    \label{eq:fish4}
\overline{b}&=&\limsup_{k\ra \infty}n_k\theta_{n_k}\leq \limsup_{\delta,\vep,\vep'\downarrow 0}\limsup_{k\ra \infty}A_{k,\delta,\vep,\vep'} + \limsup_{\delta,\vep,\vep'\downarrow 0}\limsup_{k\ra \infty}B_{k,\delta,\vep,\vep'} +\limsup_{\delta,\vep\downarrow 0}\limsup_{k\ra \infty} D_{k,\delta,\vep},
    \end{eqnarray}
where the limits are taken in the order $k\ra \infty$, $\vep'\downarrow 0$, $\vep\downarrow 0$, $\delta\downarrow 0$.

The term $A_{k,\delta,\vep,\vep'}$ can be rewritten as
\[\frac{1}{AV}\mu_{n_k}(X_{1-\delta}^{(n_k)}(1)>c\vep,X_{1}^{(n_k)}(1)>c\vep')\ra \frac{1}{AV}\N_0(X_{1-\delta}(1)>c\vep,X_{1}(1)>c\vep'), \quad \text{ as }k\ra \infty,\]
by Lemma \ref{lem:HPmass}.
  Letting $\vep'\downarrow 0$ and then $\vep\downarrow 0$ this converges to
    \[
    \frac{1}{AV}\N_0(X_{1-\delta}(1)>0,X_{1}(1)>0)=\frac{1}{AV}\N_0(X_{1}(1)>0)=2/AV,
    \]
which, in particular, does not depend on $\delta$.

Further, using Condition \ref{cond-self-repel}, the term $D_{k,\delta,\vep}$ satisfies
    \eqan{
    \label{no-resurrection}
    D_{k,\delta,\vep}&=n_k\mE\big[I_{\{0<N_{(1-\delta)n_k}\le \vep n_k\}}\mP(N_{n_k}>0|\mc{F}_{(1-\delta)n_k})\big]\leq C_{\theta}\vep n_k \theta_{\delta n_k} n_k\theta_{(1-\delta)n_k}
    \leq \frac{C\vep}{\delta(1-\delta)},\nn
    }
uniformly in $k$, since $n\theta_n$ is bounded above uniformly in $k$.
Letting $\vep\downarrow 0$, this converges to $0$.

We are left to investigate
$B_{k,\delta,\vep,\vep'}$, for which we define,
for each $m$, the measure ${\mathbb Q}_{m}=\prob(\cdot \mid N_m>0)$.  Then,
we can rewrite
    \[
    B_{k,\delta,\vep,\vep'}=n_k\theta_{(1-\delta)n_k}{\mathbb Q}_{(1-\delta)n_k}(N_{(1-\delta)n_k}>\vep n_k,0<N_{n_k}\le \vep' n_k ).
    \]
Thus, since $n_k\theta_{(1-\delta)n_k}$ is bounded above by $\frac{C}{1-\delta}\le 2C$ for $\delta<\hlf$ (where $C$ is independent of $\delta$), proving that $\limsup_{\delta,\vep,\vep'\downarrow 0}\limsup_{k\ra \infty}B_{k,\delta,\vep,\vep'}=0$ is equivalent to proving that
    \eqn{
    \label{eq:banana}
    \limsup_{\delta,\vep,\vep'\downarrow 0}\limsup_{k\ra \infty}{\mathbb Q}_{(1-\delta)n_k}(N_{(1-\delta)n_k}>\vep n_k,0<N_{n_k}\le \vep' n_k)=0.
    }

To prove \eqref{eq:banana}, we note that, for any integers $\ell_1,\ell_2\geq 0$ such that $\ell_1+\ell_2\geq 1$,
    \eqan{
    \label{mom-a}
    \expec_{{\mathbb Q}_{(1-\delta)n_{k}}}\Big[\Big(N_{(1-\delta)n_k}/n_k\Big)^{\ell_1}\Big(N_{n_k}/n_k\Big)^{\ell_2}\Big]
    &=\frac{1}{\theta_{(1-\delta)n_k}} \expec\Big[\Big(N_{(1-\delta)n_k}/n_k\Big)^{\ell_1}\Big(N_{n_k}/n_k\Big)^{\ell_2}\Big]\\
    &=\frac{1}{n_k\theta_{(1-\delta)n_k}} n_k^{-(\ell_1+\ell_2-1)}\expec[N_{(1-\delta)n_k}^{\ell_1}N_{n_k}^{\ell_2}]\nn\\
    &=\frac{1}{n_k\theta_{(1-\delta)n_k}} n_k^{-(\ell_1+\ell_2-1)} \hat{t}^{\sss(\ell_1+\ell_2+1)}_{\vec{n}_k}(0),\nn
    }
where we use that $N_{(1-\delta)n_k}>0$ when $N_{n_k}>0$, and where
$\vec{n}_k$ denotes a vector with precisely $\ell_1$ coordinates equal to $(1-\delta)n_k$ and $\ell_2$
coordinates equal to $n_k$. By Condition \ref{cond-rpoint}(a),
    \eqan{
    n_k^{-(\ell_1+\ell_2-1)} \hat{t}^{\sss(\ell_1+\ell_2+1)}_{\vec{n}_k}(0)
    &\rightarrow
    A(VA^2)^{\ell_1+\ell_2-1}\expec_{\mathbb{N}_0}\left[X_{1-\delta}(1)^{\ell_1}X_{1}(1)^{\ell_2}\right]\\
    &=\frac{2}{AV(1-\delta)} \expec_{\mathbb{N}_0}\left[\Big(VA^2X_{1-\delta}(1)\Big)^{\ell_1}
    \Big(VA^2X_{1}(1)\Big)^{\ell_2}\Big| X_{1-\delta}(1)>0\right],\nn
    }
where the last equality follows from the fact that
$\N_0(X_{1-\delta}(1)>0)=2/(1-\delta)$.
Therefore, also using that $(1-\delta)n_k\theta_{(1-\delta)n_k}\rightarrow b_{\delta}$,
    \eqan{
    \label{mom-b}
    \expec_{{\mathbb Q}_{(1-\delta)n_{k}}}\Big[\Big(N_{(1-\delta)n_k}/n_k\Big)^{\ell_1}\Big(N_{n_k}/n_k\Big)^{\ell_2}\Big]
    &\rightarrow \frac{2}{AV b_{\delta}} \expec_{\mathbb{N}_0}\left[\Big(VA^2 X_{1-\delta}(1)\Big)^{\ell_1}
    \Big(VA^2X_{1}(1)\Big)^{\ell_2}\Big| X_{1-\delta}(1)>0\right].
    }
We recognize the above joint moments as the joint moments of $(X,Y)$ with distribution $(1-\alpha_{\delta})\delta_{(0,0)}+\alpha_{\delta} \nu_{\delta}$, where $\delta_{(0,0)}$
is the point measure on the vector $(0,0)$ and $\nu_{\delta}$ is the
law of $(A^2VX_{1-\delta}(1),A^2VX_{1}(1))$ under $\N_0(\cdot | X_{1-\delta}(1)>0)$,
and with $\alpha_{\delta}=2/(AV b_{\delta})$.
For any $t>1-\delta$,
    \eqn{
    \N_0(X_{t}(1)=0|X_{1-\delta}(1)>0)=1-(1-\delta)/t,
    }
so that
    \eqn{
    \label{X1=0}
    \nu_{\delta}(X_{1}(1)=0)=1-(1-\delta)=\delta.
    }
Let $(X_n, Y_n)$ be a two-dimensional distribution. Again by
\cite[Theorems 30.1 and 30.2, and Problems 30.5 and 30.6]{Bill86},
convergence of the joint moments of $(X_n,Y_n)$ to those of
$(X,Y)$ implies convergence in distribution when the
moment generating function of both $X$ and $Y$ are finite in
a neighborhood of 0.  Under the conditional law
$\N_0(\cdot | X_{1-\delta}(1)>0)$, the distribution
of $A^2VX_{1-\delta}(1)$ is exponential with mean $(1-\delta)A^2V/2$
(see e.g., \cite[Theorem 1.4]{HHS_IIC}), and by \eqref{X1=0}, $A^2VX_1(1)$
is 0 with probability $\delta$ and an exponential with mean $A^2V/2$
with probability $1-\delta$.
As a result, the distribution of both limits $X$ and $Y$ are mixtures of point masses
at 0 with probabilities $1-\alpha_{\delta}$ and $1-\alpha_{\delta}+\alpha_{\delta}\delta$
and exponentials with positive means $\lambda_{\sss X}$ and $\lambda_{\sss Y}$. Therefore,
their moment generating functions are finite in a neighborhood of zero, so that
$\big(N_{(1-\delta)n_k}/n_k, N_{n_k}/n_k\big)$ converges in distribution
to $(X,Y)$ having distribution $(1-\alpha_{\delta})\delta_{(0,0)}+\alpha_{\delta} \nu_{\delta}$.

Thus, as $k\ra \infty$,
    \eqan{
    {\mathbb Q}_{(1-\delta)n_{k}}(N_{(1-\delta)n_{k}}>\vep n_{k}, N_{n_k}\leq \vep'n_k)
    &\rightarrow
    \alpha_{\delta}\nu_{\delta}(A^2V X_{1-\delta}(1)>\vep , A^2VX_{1}(1)\leq \vep').\nn
    }
When $\vep'\downarrow 0$,
    \eqn{
    \nu_{\delta}(A^2V X_{1-\delta}(1)>\vep , A^2VX_{1}(1)\leq \vep')\rightarrow
    \nu_{\delta}(X_{1-\delta}(1)>\vep c, X_{1}(1)=0)
    \leq
    \nu_{\delta}(X_{1}(1)=0)
    =\delta,
    }
where we use \eqref{X1=0}.  Letting $\delta\downarrow 0$, we obtain
\eqref{eq:banana}. We conclude that $\limsup_{n\rightarrow \infty} n\theta_n=\bar{b}\leq 2/(AV)$,
which, together with \eqref{weak-lb}, shows that $\lim_{n\rightarrow \infty} n\theta_n=2/(AV)$,
as required.


The fact that, conditionally on $N_{nt}>0$,
the finite-dimensional distributions of
$(N_{sn}/n)_{s\geq t}$ converge to those
of Feller's branching diffusion started from an exponential
random variable with mean $A^2Vt/2$ can be obtained as follows.  Fix $t=s_0<s_1<\dots <s_r<\infty$, and let $0<\ell=\sum_{j=0}^r\ell_j$, where $\ell_j$, $j=0,\dots, r$ are non-negative integers.  Set $\vec{s}=(s_0,\dots, s_r)$.  As in \eqref{mom-a} we have
 \eqan{
    \expec_{{\mathbb Q}_{tn}}\Bigg[\prod_{j=0}^r\Big(\frac{N_{s_jn}}{n}\Big)^{\ell_j}\Bigg]
    &=\frac{1}{n\theta_{tn}}\, n^{1-\ell}\, \hat{t}^{\sss(\ell+1)}_{\vec{s}n}(0),\nn
    }
where we now know that $tn\theta_{tn}\rightarrow 2/(AV)$, and as before
\eqan{
    n^{1-\ell}\, \hat{t}^{\sss(\ell+1)}_{\vec{s}n}(0)
    &\rightarrow  \frac{2}{AVt} \expec_{\mathbb{N}_0}\left[\prod_{j=0}^r\Big(VA^2 X_{s_j}(1)\Big)^{\ell_j}\Big| X_{t}(1)>0\right].\nn
    }
Thus the joint moments converge as in \eqref{mom-b}, i.e.
    \eqan{
    \label{mom-b*}
    \expec_{{\mathbb Q}_{tn}}\Bigg[\prod_{j=0}^r\Big(\frac{N_{s_jn}}{n}\Big)^{\ell_j}\Bigg]
    &\rightarrow \expec_{\mathbb{N}_0}\left[\prod_{j=0}^r\Big(VA^2 X_{s_j}(1)\Big)^{\ell_j}\Big| X_{t}(1)>0\right].}
Finally, the fact that $(A^2VX_{s}(1))_{s\geq t}$ is Feller's branching diffusion
follows from \cite{Fell51}. Again by
\cite[Theorems 30.1 and 30.2, and Problems 30.5 and 30.6]{Bill86},
and the above bound on the moment generating function of $X_{s}(1)$,
this completes the proof.
\qed

\paragraph{Acknowledgements.}
The work of RvdH was supported
in part by the Netherlands Organisation
for Scientific Research (NWO).
We thank Akira Sakai for discovering an error in a previous
version, and Tim Hulshof for useful suggestions that helped
us to improve the presentation.


\bibliographystyle{plain}
\def\cprime{$'$}

\end{document}